\documentclass[aap,MSNbibl,dvips]{arximspdf}

\doi{10.1214/12-AAP896} 
\volume{23}
\issue{5}
\pubyear{2013}
\firstpage{2139}
\lastpage{2160}

\makeatletter

\newcommand{\mmbox}[1]{\mbox{\fontsize{8.36pt}{10pt}\selectfont{\textup{#1}}}}

\newcommand{\rrvert}{\vert}
\newcommand{\llvert}{\vert}

\newtheorem{theorem}{Theorem}
\newproclaim{example}[theorem]{Example}
\newtheorem{lemma}[theorem]{Lemma}
\newproclaim{remark}[theorem]{Remark}

\newcommand{\R}{\mathbb{R}}
\newcommand{\N}{\mathbb{N}}
\newcommand{\E}{\mathbb{E}}
\newcommand{\F}{\mathcal{F}}
\renewcommand{\P}{\mathbb{P}}

\newcommand{\rp}{\eta}
\newcommand{\rrp}{{\bolds{\eta}}}

\newcommand{\Lip}{\operatorname{Lip}}
\newcommand{\mS}{\mathcal{S}}

\newcommand{\mY}{\mathcal{Y}}

\newcommand{\olphi}{\phi}
\newcommand{\olpsi}{\psi}

\newcommand{\tS}{\tilde{S}}
\newcommand{\x}{\mathbf{x}}

\newcommand{\Holl}{\mbox{\fontsize{8.36pt}{10pt}\selectfont{H\"{o}l}}}
\newcommand{\Holll}{\mbox{\fontsize{6.6pt}{8pt}\selectfont{\textup{H\"{o}l}}}}

\makeatother

\begin{document}
\begin{frontmatter}

\title{Robust filtering: Correlated noise and multidimensional observation}
\runtitle{Robust filtering}

\begin{aug}
\author[A]{\fnms{D.} \snm{Crisan}\thanksref{t1}\ead[label=e1]{d.crisan@ic.ac.uk}},
\author[B]{\fnms{J.} \snm{Diehl}\thanksref{t2}\ead[label=e2]{diehl@math.tu-berlin.de}},
\author[C]{\fnms{P. K.} \snm{Friz}\corref{}\thanksref{t3}\ead[label=e4]{friz@math.tu-berlin.de}\ead[label=e5]{friz@wias-berlin.de}}
\and
\author[B]{\fnms{H.} \snm{Oberhauser}\thanksref{t3}\ead[label=e3]{h.oberhauser@gmail.com}}
\runauthor{Crisan, Diehl, Friz and Oberhauser}
\affiliation{Imperial College London, Technical University Berlin, Technical
University Berlin and
Weierstrass Institute for Applied Analysis and Stochastics Berlin, and~Technical University Berlin}
\address[A]{D. Crisan\\
Department of Mathematics\\
Imperial College London\\
180 Queen's Gate\\
London SW7 2BZ\\
United Kingdom\\
\printead{e1}} 
\address[B]{J. Diehl\\
H. Oberhauser\\
Institut f\"{u}r Mathematik, MA 7-2\\
Technical University Berlin\\
Strasse des 17. Juni 136\\
10623 Berlin\\
Germany\\
\printead{e2}\\
\hphantom{E-mail: }\printead*{e3}}
\address[C]{P. Friz\\
Institut f\"{u}r Mathematik, MA 7-2\\
Strasse des 17. Juni 136\\
10623 Berlin\\
Germany\\
and\\
Weierstrass-Institut f\"{u}r Angewandte Analysis\\
\quad und Stochastik\\
Mohrenstrasse 39\\
10117 Berlin\\
Germany\\
\printead{e4}\\
\hphantom{E-mail: }\printead*{e5}}
\end{aug}

\thankstext{t1}{Supported in part by the
EPSRC Grant EP/H0005500/1.}

\thankstext{t2}{Supported by DFG Grant SPP-1324.}

\thankstext{t3}{Supported by the European Research Council
under the European Union's Seventh Framework Programme
(FP7/2007-2013)/ERC Grant agreement 258237.}

\received{\smonth{1} \syear{2012}}
\revised{\smonth{8} \syear{2012}}

%
\begin{abstract}
In the late seventies, Clark [In \textit{Communication Systems and
Random Process Theory} (\textit{Proc. 2nd NATO Advanced Study Inst.},
\textit{Darlington}, 1977) (1978) 721--734, Sijthoff \& Noordhoff]
pointed out that it would be natural for $\pi_t$, the solution of the
stochastic filtering problem, to depend continuously on the observed
data $Y=\{Y_s,s\in[0,t]\}$. Indeed, if the signal and the observation
noise are independent one can show that, for any suitably chosen test
function $f$, there exists a continuous map $\theta^f_t$, defined on
the space of continuous paths $C([0,t],\R^d)$ endowed with the uniform
convergence topology such that $\pi_t(f)=\theta^f_{t}(Y)$, almost
surely; see, for example, Clark [In \textit{Communication Systems and
Random Process Theory} (\textit{Proc. 2nd NATO Advanced Study Inst.},
\textit{Darlington}, 1977) (1978) 721--734, Sijthoff \& Noordhoff],
Clark and Crisan [\textit{Probab. Theory Related Fields} \textbf{133}
(2005) 43--56], Davis [\textit{Z.~Wahrsch. Verw. Gebiete} \textbf{54}
(1980) 125--139], Davis [\textit{Teor. Veroyatn. Primen.} \textbf{27}
(1982) 160--167], Kushner [\textit{Stochastics} \textbf{3} (1979)
75--83]. As shown by Davis and Spathopoulos [\textit{SIAM J. Control
Optim.} \textbf{25} (1987) 260--278], Davis [In \textit{Stochastic
Systems}: \textit{The Mathematics of Filtering and Identification and
Applications}, \textit{Proc. NATO Adv. Study Inst. Les Arcs, Savoie,
France 1980} 505--528], [In \textit{The Oxford Handbook of Nonlinear
Filtering} (2011) 403--424 Oxford Univ. Press], this type of
\textit{robust} representation is also possible when the signal and the
observation noise are correlated, provided the observation process is
scalar. For a general correlated noise and multidimensional
observations such a representation does not exist. By using the theory
of rough paths we provide a solution to this deficiency: the
observation process $Y$ is ``lifted'' to the process $\mathbf{Y}$ that
consists of $Y$ and its corresponding L\'evy area process, and we show
that there exists a continuous map $\theta_{t}^f$, defined on a
suitably chosen space of H\"older continuous paths such that
$\pi_t(f)=\theta_{t}^f(\mathbf{Y})$, almost surely.
\end{abstract}

%

\end{frontmatter}

\section{Introduction}\label{sec1}

Let $(\Omega, \F, (\F_t)_{t\ge0}, \P)$ be a filtered probability
space on which we have defined a two-component diffusion process $(X,Y)$
solving a stochastic differential equation driven by a multidimensional
Brownian motion.
One assumes that the first component $X$ is unobservable, and the
second component $Y$ is observed. The filtering problem consists of
computing the conditional distribution of the unobserved component, called
the \textit{signal} process, given the \textit{observation} process $Y$.
Equivalently, one is interested in computing
\[
\pi_t(f)=\mathbb{E}\bigl[f(X_t,Y_t)|
\mathcal{Y}_t\bigr],
\]
where $\mathcal{Y}=\{\mathcal{Y}_t,t\ge0\}$ is the observation
filtration, and $f$ is a suitably chosen test function.
An elementary measure theoretic result tells us\setcounter{footnote}{3}\footnote{See, for
example, Proposition 4.9, page 69, in~\cite{breiman}.} that there exists
a Borel-measurable map $\theta^f_{t}\dvtx C([0,t],\R^{d_Y})\to\R$, such
that
%
\begin{equation}\label{nonrobust}
\pi_t( f )=\theta^f_{t}(Y_\cdot),\qquad
\mathbb{P}\mbox{-a.s.,}
\end{equation}
where $d_Y$ is the dimension of the observation state space, and
$Y_{\cdot}$ is the path-valued random variable
\[
Y_\cdot\dvtx\Omega\rightarrow C\bigl([0,t],\R^{d_Y}\bigr),\qquad
Y_\cdot(\omega)=\bigl(Y_s(\omega), 0\leq s\leq t\bigr).
\]
Of course, $\theta^f_{t}$ is not unique. Any other function $\bar
\theta^f_{t}$
such that
\[
\mathbb{P}\circ Y_\cdot^{-1} \bigl( \bar\theta^f_{t}
\neq\theta^f_{t} \bigr) =0,
\]
where $\mathbb{P}\circ Y_\cdot^{-1}$ is the distribution of $Y_{\cdot
}$ on
the path space $C([0,t],\R^{d_Y})$ can replace $\theta^f_{t}$ in
(\ref
{nonrobust}). It would be desirable to solve this ambiguity by choosing
a suitable representative from the class of functions that satisfy
(\ref
{nonrobust}). A \textit{continuous} version, if it exists, would enjoy
the following uniqueness property: if the law of the observation
$\mathbb{P}\circ Y_\cdot^{-1}$ positively charges all nonempty
open sets in $C([0,t],\R^{d_Y})$, then there exists a
unique continuous function $\theta^f_{t}$ that satisfies (\ref
{nonrobust}). In this case, we call $\theta^f_{t}(Y_\cdot)$ the
\textit{robust version} of $\pi_t( f )$ and equation (\ref{nonrobust}) is
the robust representation formula for the solution of the stochastic
filtering problem.

The need for this type of representation arises when the filtering framework
is used to model and solve ``real-life'' problems. As explained in a
substantial number of papers (e.g.,
\cite{clark,cc,d1,d2,d3,d5,d4,k})
the model chosen for the ``real life'' observation process $\bar{Y}$
may not be a perfect
one. However,\vadjust{\goodbreak}
if $\theta^f$ is continuous (or even locally Lipschitz, as in the
setting of~\cite{cc}),
and
as long as the distribution of $\bar{Y}_{\cdot}$ is close
in a weak sense to that of $Y_\cdot$ (and some integrability
assumptions hold), the estimate $\theta^f_{t}( \bar{Y}_\cdot) $
computed on the actual observation will still be reasonable, as
$\mathbb{E}[(f(X_t,Y_t)-\theta^f_{t}(\bar{Y}_\cdot))^2]$ is close
to the
idealized error $\mathbb{E}[(f(X_t,Y_t)-\theta^f_{t}(Y_\cdot))^2]$.

Moreover, even when $Y$ and $\bar{Y}$ actually coincide, one is never
able to obtain and
exploit a continuous stream of data as modeled by the continuous
path $Y_\cdot(\omega)$. Instead the observation arrives and
is processed at discrete moments in time
\[
0=t_0<t_1<t_2<\cdots<t_n=t.
\]
However, the continuous path $\hat{Y}_\cdot(\omega)$
obtained from the discrete observations $(Y_{t_i}(\omega))_{i=1}^{n}$
by linear interpolation is close to $Y_\cdot(\omega)$
(with respect to the supremum norm on $C([0,t],\R^{d_Y}) $);
hence, by the same argument, $\theta^f_{t}(\hat{Y}_\cdot)$
will be a sensible approximation to $\pi_t( f )$. To conclude the
discussion on the un-correlated framework, let us also mention that
Kushner introduces in~\cite{kushner1980robust} a robust computable
approximation for the filtering
solution.

In the following, we will assume that the pair of processes $(X,Y)$
satisfy the equation
%
\begin{eqnarray}
\label{eqMainSDEIto'}
dX_t &=& l_0(X_t,Y_t)
\,dt + \sum_k Z_k( X_t,Y_t
) \,dW^k_t + \sum_j
L_j(X_t,Y_t) \,dB^j_t,
\\
\label{eqMainSDEItoY'}
dY_t &=& h(X_t,Y_t)\,dt + dW_t
\end{eqnarray}
with $X_0$ being a bounded random variable and $Y_0 = 0$. In (\ref
{eqMainSDEIto'}) and (\ref{eqMainSDEItoY'}), the process $X$ is the
$d_X$-dimensional signal, $Y$ is the $d_Y$-dimensional observation, $B$
and $W$ are independent $d_B$-dimensional, respectively,
$d_Y$-dimensional Brownian motions independent of $X_0$. Suitable
assumptions on the coefficients
$l_0, L_1,\ldots, L_{d_B}\dvtx \R^{d_X+d_Y} \to\R^{d_X}$,
$Z_1,\ldots, Z_{d_Y}\dvtx \R^{d_X+d_Y} \to\R^{d_X}$ and
$h = (h^1,\ldots, h^{d_Y})\dvtx \R^{d_X+d_Y} \to\R^{d_Y}$ will be
introduced later on.
This framework covers a wide variety of applications of stochastic
filtering (see, e.g.,~\cite{cr} and the references therein) and has
the added advantage that, within it, $\pi_t( f )$ admits an
alternative representation that is crucial for the construction of its
robust version. Let us detail this representation first.

Let $u=\{u_{t},t>0\}$ be the process defined by
%
\begin{equation}\label{ut}
u_t=\exp\Biggl[ - \sum_{i=1}^{d_Y}
\biggl( \int_0^t h^i(X_s,Y_s)
\,dW^{i}_s - \frac{1}{2} \int_0^t
\bigl(h^i(X_s,Y_s)\bigr)^2 \,ds
\biggr) \Biggr].
\end{equation}
Then, under suitable assumptions,\footnote{For example, if Novikov's
condition is satisfied, that is, if
$ \mathbb{E} [ \exp( \frac{1}{2}\int_{0}^{t}\|
h^{i}(X_s,\break Y_s)\|^{2}\,ds
) ] <\infty\label{Novikov}$ for all $t>0$, then $u$ is
a martingale.
In particular it will be satisfied in our setting, in which $h$ is bounded.}
$u$ is a martingale which is used to construct the probability measure
$\P_0$ equivalent to $\mathbb{P}$ on
$\bigcup_{0 \le t < \infty} \mathcal{F}_t$ whose Radon--Nikodym
derivative with respect to $\mathbb{P}$
is given by $u$, namely,
\[
\frac{d \mathbb{P}_0}{d \mathbb
{P}}\bigg|_{\mathcal{F}_t}=u_t.
\]
Under $\P_0$, $Y$ is a Brownian motion independent of $B$.
Moreover the equation for the signal process $X$ becomes
%
\begin{equation}
dX_t \label{eqMainSDEIto''} =
\bar{l}_0(X_t,Y_t) \,dt + \sum
_k Z_k( X_t,Y_t )
\,dY^k_t + \sum_j
L_j(X_t,Y_t) \,dB^j_t.
\end{equation}
Observe that equation (\ref{eqMainSDEIto''}) is now written in terms
of the pair of Brownian motions $(Y,B)$ and the coefficient $\bar
{l}_0$ is given by $\bar{l}_0=l_0+\sum_k Z_kh_k$. Moreover, for any
measurable, bounded function $f\dvtx\R^{d_X+d_Y} \to\R$, we have the
following formula, called the Kallianpur--Striebel formula:
%
\begin{equation}
\label{eqPtf''} \pi_t(f) =
\frac{p_t(f)}{p_t(1)},\qquad p_t(f):= \E_0\bigl[
f(X_t,Y_t) v_t| \mathcal{Y}_t
\bigr], 
\end{equation}
where $v=\{v_{t},t>0\}$ is the process defined as $v_t:=\exp( I_t )$,
$t\ge0$ and
%
\begin{equation}
\label{eqPtf'} I_t:= \sum
_{i=1}^{d_Y} \biggl( \int_0^t
h^i(X_r,Y_r) \,dY^i_r
- \frac{1}{2} \int_0^t
\bigl(h^i(X_r,Y_r)\bigr)^2 \,dr
\biggr),\qquad t\ge0.
\end{equation}
The representation (\ref{eqPtf''}) suggests the following three-step
methodology to construct a robust representation formula for $\pi_t^f$:
\begin{longlist}[\textit{Step} 2.]
\item[\textit{Step} 1.] We construct the triplet of processes
$(X^y,Y^y,I^y)$\footnote{As we shall see momentarily, in the
uncorrelated case the choice of $Y^y$ will trivially be $y$.
In the correlated case we make it part of the \textit{SDE with rough
drift}, for (notational) convenience.}
corresponding to the pair $(y,B)$ where $y$ is now a \textit{fixed}
observation path $y_.=\{y_s,s\in[0,t]\}$ belonging to a suitable class
of continuous functions and prove that the random variable $f(X^y, Y^y)
\exp( I^y )$ is $\P_0$-integrable.
\item[\textit{Step} 2.]
We prove that the function $y_{\cdot}\to g^f_t(y_{\cdot})$ defined as
%
\begin{equation}
\label{genfor'} g^f_t(y_.)=
\E_0 \bigl[ f\bigl(X^y_t,Y^y_t
\bigr) \exp\bigl( I^y_t \bigr) \bigr]
\end{equation}
is continuous.
\item[\textit{Step} 3.]
We prove that $g^f_t(Y_{\cdot})$ is a version of $p_t(f)$. Then, following
(\ref{eqPtf''}), the function, $y_{\cdot}\to\theta^f_t(y_{\cdot})$ defined as
%
\begin{equation}
\label{genfor}
\theta^f_{t}={g^f_t\over g^1_t}
\end{equation}
provides the robust version of of $\pi_t(f)$.
\end{longlist}
We emphasize that step 3 cannot be omitted from the
methodology. Indeed one has to prove that $g^f_t(Y_{\cdot})$ is a version of
$p_t(f)$ as this fact is not immediate from the definition of $g^f_t$.

\textit{Step} 1 is immediate in the particular case when only the
Brownian motion $B$ drives $X$ (i.e., the coefficient $Z=0$) and $X$ is
itself a diffusion, that is, it satisfies an equation of the form
%
\begin{equation}
dX_t \label{simple} = l_0(X_t) \,dt + \sum
_j L_j(X_t)
\,dB^j_t,
\end{equation}
and $h$ does only depend on $X$.
In this case the process $(X^y, Y^y)$ can be taken to be the pair
$(X,y)$. Moreover, we can define $I^y$ by the formula
%
\begin{equation}
\label{eqPtf'''} I^y_t
:= \sum_{i=1}^{d_Y} \biggl(h^i(X_t)y^i_t-
\int_0^t y^i_r\,dh^i(X_r)-
\frac
{1}{2} \int_0^t
\bigl(h^i(X_r,Y_r)\bigr)^2 \,dr
\biggr),\qquad t\ge0,\hspace*{-28pt}
\end{equation}
provided the processes $h^i(X)$ are semi-martingales. In (\ref
{eqPtf'''}), the integral $\int_0^t y^i_r\,dh^i(X_r)$ is the It\^o
integral of the nonrandom process $y^{i}$ with respect to $h^i(X)$.
Note that the formula for
$I_t^{y}$ is obtained by applying integration by parts to the
stochastic integral in (\ref{eqPtf'})
%
\begin{equation}
\label{intpart} \int_0^t h^i(X_r)
\,dY^i_r=h^i(X_t)Y^i_r-
\int_0^t Y^i_r\,dh^i(X_r),
\end{equation}
and replacing the process $Y$ by the fixed path $y$ in (\ref
{intpart}). This approach has been successfully used to study the
robustness property for the filtering problem for the above case in a
number of papers~\cite{clark,cc,k}.

The construction of the process $(X^y,Y^y,I^y)$ is no longer immediate
in the case when $Z\ne0$, that is, when the signal is driven by both
$B$ and $W$ (the correlated noise case). In the case when the
observation is one-dimensional, one can solve this problem by using a
method akin with the Doss--Sussmann ``pathwise solution'' of a
stochastic differential equation; see~\cite{doss,sussmann}.
This approach has been employed by Davis to extend the robustness
result to the correlated noise case with scalar observation; see \cite
{d1,d2,d3,d4}. In this case one constructs first a diffeomorphism which
is a pathwise solution of the equation\footnote{Here $d^Y=1$ and $Y$ is scalar.}
%
\begin{equation}
\label{difeo} \phi(t,x)=x + \int_0^t Z\bigl(
\phi(s,x) \bigr)\circ dY_t.
\end{equation}
The diffeomorphism is used to express the solution $X$ of equation
(\ref{eqMainSDEIto''}) as a composition between the diffeomorphism
$\phi$ and the solution of a stochastic differential equation driven\vadjust{\goodbreak}
by $B$ only and whose coefficients depend continously on $Y$. As a
result, we can make sense
of $X^y$.
$I^y$ is then defined by a suitable (formal) integration by parts that
produces a pathwise interpretation of the stochastic integral appearing
in (\ref{eqPtf'}), and $Y^y$ is chosen to be $y$, as before.
The robust representation formula is then introduced as per (\ref{genfor}).
Additional results for the correlated noise case with scalar
observation can be found in~\cite{florchinger1993zakai}.
The extension of the robustness result to special cases of the
correlated noise and multidimensional observation
has been tackled in several works.
Robustness results in the correlated setting have been obtained by
Davis in~\cite{d1,d4} and Elliott and Kohlmann in~\cite{elliott1981robust},
under a commutativity condition on the signal vector fields.
Florchinger and Nappo~\cite{florchingernappo} do not have correlated noise,
but allow the coefficients to depend on the signal and the
observation.\footnote{We thank the anonymous referee for these references.}
To sum up, all previous works on the robust representation problem
either treat the uncorrelated case, the case with one-dimensional
observation or the case where the Lie brackets of the relevant vector
fields vanish.
In parallel, Bagchi and Karandikar treat in~\cite{bagchi1994white}
a different model with ``finitely additive'' state white noise and
``finitely additive'' observation noise.
Robustness there is virtually built into the problem.


An alternative framework is that where the signal and the observation
run in discrete time. In this case the filtering problem is well
understood and has been studied in many works, including the monograph
\cite{del2004feynman} and the articles
\cite{del2001monte,del2010concentration,del2000branching}. These works
include an analysis of discrete time filtering problems and their
approximation models, including particle approximation, approximate
Bayesian computation, filtering models, etc. We note that in this
context the continuity of the filter with respect to the observation
data holds true\footnote{which can be easily seen using the
representation of Lemma 2.1 in~\cite{del2000branching}.} provided very
natural conditions are imposed on the model: for example, the
likelihood functions are assumed to be continuous and bounded (which
includes the Gaussian case).

To our knowledge, the general correlated noise and multidimensional
observation case has not been studied, and it is the subject of the
current work. In this case it turns out that we cannot hope to have
robustness in the sense advocated by Clark. More precisely, there may
not exist a map continuous map
$\theta^f_t\dvtx C([0,t],\R^{d_Y})\to\R$, such that the representation
(\ref{nonrobust}) holds almost surely. The following is a simple
example that illustrates this.
%
\begin{example}
Consider the filtering problem where the signal and the observation
process solve the following pair of equations:
\begin{eqnarray*}
X_t &=& X_0 + \int_0^t
X_r \,d \bigl[ Y^1_r + Y^2_r
\bigr] + \int_0^t X_r \,dr,
\\
Y_t &=& \int_0^t
h(X_r) \,dr + W_t,
\end{eqnarray*}
where $Y$ is two-dimensional and $\P(X_0=0) = \P(X_0=1) = \frac{1}{2}$.
Then with $f,h$ such that $f(0) = h_1(0)=h_2(0) = 0$ one can explicitly compute
%
\begin{eqnarray}
\label{example}\qquad
&&\E\bigl[ f( X_t ) | \mY_t \bigr] \nonumber\\[-8pt]\\[-8pt]
&&\qquad=
{f( \exp( Y^1_t + Y^2_t ) ) \over1 + \exp( -
\sum_{k=1,2} \int_0^t h^k( \exp( Y_r ) ) \,dY^k_r
+\int_0^t \|h(\exp(Y_r))\|^2 \,dr/2 )}.\nonumber
\end{eqnarray}
Following the findings of rough path theory (see, e.g., \cite
{bibFrizVictoir,bibLyons94,bibLyonsCaruanaLevy04,bibLyonsQian02}) the
expression on the right-hand side of (\ref{example}) is not continuous
in supremum norm
(nor in any other metric on path space)
because of the stochastic integral. Explicitly, this follows, for
example, from Theorem 1.1.1 in~\cite{bibLyonsQian02} by rewriting the
exponential term as the solution to a stochastic differential equation
driven by $Y$.
\end{example}
Nevertheless, we can show that a variation of the robustness
representation formula still exists in this case. For this we need to
``enhance'' the original process $Y$ by adding a second component to it
which consists of its iterated integrals (that, knowing the path, is in
a one-to-one correspondance with the L{\'e}vy area process).
Explicitly we consider the process $\mathbf{Y}=\{\mathbf
{Y}_t,t\ge0\}$ defined as
%
\begin{equation}
\label{ybold} \mathbf{Y}_t = \left(Y_t, \pmatrix{
\displaystyle \int_0^t Y^1_r \circ
dY^1_r & \cdots& \displaystyle \int_0^t
Y^1_r \circ dY^{d^Y}_r
\cr
\vdots&
\vdots& \vdots
\cr
\displaystyle \int_0^t
Y^{d^Y}_r \circ dY^1_r & \cdots&
\displaystyle \int_0^t Y^{d^Y}_r \circ
dY^{d^Y}_r} \right),\qquad t\ge0.
\end{equation}
The stochastic integrals in (\ref{ybold}) are Stratonovich integrals.
The state space of $\mathbf{Y}$ is
$G^{2}(\R^{d_Y}) \cong\R^{d_Y} \oplus\operatorname{so}(d_Y)$,
where $\operatorname{so}(d_Y)$ is the set of anti-symmetric matrices
of dimension
$d_Y$.\footnote{More generally, $G^{[1/\alpha]}(\R^d)$ is the
``correct'' state space for a geometric $\alpha$-H\"older rough path;
the space of such paths subject to $\alpha$-H\"older regularity (in
rough path sense) yields a complete metric space under $\alpha$-H\"
older rough path metric. Technical details of geometric rough path
spaces (as found, e.g., in Section 9 of~\cite{bibFrizVictoir}) are not
required for understanding the results of the present paper.}
Over this state space we consider not the space of continuous function,
but a subspace $\mathcal{C}^{0,\alpha}$ that contains paths $\eta
\dvtx[0,t]\rightarrow
G^{2}( \R^{d_Y} )$ that are $\alpha$-H\"{o}lder in the $\R
^{d_Y}$-component and somewhat ``$2\alpha$-H\"{o}lder'' in the $
\operatorname{so}(d_Y)$-component, where $\alpha$ is a suitably
chosen constant $\alpha<1/2$. Note that there exists a modification of
$\mathbf{Y}$ such that $\mathbf{Y}(\omega)\in\mathcal
{C}^{0,\alpha}$ for all
$\omega$ (Corollary 13.14 in~\cite{bibFrizVictoir}).

The space $\mathcal{C}^{0,\alpha}$ is endowed with the $\alpha$-H\"
older rough path
metric under which
$\mathcal{C}^{0,\alpha}$ becomes a complete metric space. The
\textit{main result of
the paper} (captured in Theorems~\ref{thmContinuous} and \ref
{thmVersionOfConditionalExpectation}) is that there exists a continuous
map $\theta_t^f\dvtx\mathcal{C}^{0,\alpha}\to\R$, such that
%
\begin{equation}\label{newrobust}
\pi_t( f )=\theta_t^f(\mathbf{Y}_\cdot),\qquad
\mathbb{P}\mbox{-a.s.}
\end{equation}
Even though the map is defined on a slightly more abstract space,
it nonetheless enjoys the desirable properties described above for the
case of a continuous version on $C([0,t], \R^d)$.
Since $\P\circ\mathbf{Y}^{-1}$ positively charges all nonempty
open sets of $\mathcal{C}^{0,\alpha}$,\footnote{This fact is a
consequence of the
support theorem of Brownian motion
in H\"older rough path topology~\cite{bibFriz05}; see also Chapter 13
in~\cite{bibFrizVictoir}.}
the continuous version we construct will be unique. Also, it provides a
certain model robustness, in the sense that
$\mathbb{E}[(f(X_t)-\theta^f_{t}(\bar{\mathbf{Y}}_\cdot))^2]$ is
well approximated by the
idealized error $\mathbb{E}[(f(X_t)-\theta^f_{t}(\mathbf{Y}_\cdot))^2]$,
if $\bar{\mathbf{Y}}_\cdot$ is close in distribution to $\mathbf
{Y}_\cdot$.
The problem of discrete observation is a little more delicate.
One one hand,
it is true that the rough path lift $\hat{\mathbf{Y}}$ calculated
from the linearly interpolated Brownian motion $\hat Y$
will converge to the true rough path $\mathbf{Y}$
in probability as the mesh goes to zero (Corollary 13.21 in~\cite
{bibFrizVictoir}),
which implies that $\theta_t^f(\hat{\mathbf{Y}})$
is close in probability to $\theta_t^f( \mathbf{Y})$
(we provide local Lipschitz estimates for $\theta^f$).
Actually, most sensible approximations will do,
as is, for example, shown in Chapter 13 in~\cite{bibFrizVictoir}
(although, contrary to the uncorrelated case, not all interpolations
that converge
in uniform topology will work; see, e.g., Theorem 13.24 ibid).
But these are probabilistic statements, that somehow miss the pathwise
stability that one wants to provide with $\theta_t^f$.
If, on the other hand, one is able to observe at discrete time points
not only the process itself,
but also its second level, that is, the area,
one can construct an interpolating rough path using geodesics
(see, e.g., Chapter~13.3.1 in~\cite{bibFrizVictoir}) which is close to
the true (lifted) observation path $\mathbf{Y}$ in the relevant metric
\textit{for all realizations $\mathbf{Y} \in\mathcal{C}^{0,\alpha}$}.

The following is the outline of the paper:
In the next section, we enumerate the common notation used throughout
the paper.
In Section~\ref{sec3} we introduce the notion of a stochastic differential
equation with rough drift, which is necessary for our main result
and correspond to step 1 above.
We present it separately of the filtering problem, since we believe
this notion to be of independent interest.
The proof of the existence of a solution of a stochastic differential
equation with rough drifts and its properties is postponed to Section~\ref{sec5}.
Section~\ref{sec4} contains the main results of the paper and the assumptions
under which they hold true.
Steps 2 and 3 of above mentioned methodology are carried out in Theorems
\ref{thmContinuous} and~\ref{thmVersionOfConditionalExpectation}.

\section{Nomenclature}\label{sec2}

$\Lip^{\gamma}$ is the set of $\gamma$-Lipschitz\footnote{In the sense
of E. Stein, that is, bounded $k$th dervative for
$k=0,\ldots,\lfloor\gamma\rfloor$ and $\gamma-
\lfloor\gamma\rfloor$-H\"older continuous $\lfloor \gamma\rfloor$th
derivative.} functions $a\dvtx\R^m\to\R^n$ where $m$ and $n$ are
chosen according to the context.

$G^{2}(\R^{d_Y}) \cong\R^d \oplus\operatorname{so}(d_Y)$ is
the state space for a $d_Y$-dimensional Brownian motion (or, in general
for an arbitrary semi-martingale) and its corresponding L\'evy area.

$\mathcal{C}^{0,\alpha}:=C_0^{0, \alpha\mbox{-}\Holl}( [0,t], G^{2}(
\R^{d_Y} ) )$ is the set of
geometric $\alpha$-H\"{o}lder rough paths $\eta\dvtx[0,t]\rightarrow
G^{2}( \R^{d_Y} )$ starting at $0$.
We shall use the nonhomogenous metric $\rho_{\alpha\mbox{-}\Holl}$ on this space.

In the following we will make use of an auxiliary filtered probability
space $( \bar{\Omega}, \bar{\F}, (\bar{F}_t)_{t\ge0}, \bar{\P} )$
carrying a $d_B$-dimensional Brownian motion $\bar{B}$.\footnote{We
introduce this auxiliary probability space, since in the proof of
Theorem~\ref{thmVersionOfConditionalExpectation} it will be easier to
work on a product space separating the randomness coming from $Y$ and
$B$.
A similar approach was followed in the proof of Theorem 1 in \cite
{bibBainCrisan}.}

Let $\mS^0 = \mS^0(\bar{\Omega})$ denote the space of adapted,
continuous processes in $\R^{d_S}$,
with the topology of uniform convergence in probability.

For $q \ge1$ we denote by $\mS^q = \mS^q(\bar{\Omega})$ the space
of processes $X \in\mS^0$ such that
\[
\|X\|_{\mS^q}:= \Bigl( \bar{\E}\Bigl[ \sup_{s\le t}
|X_t|^q \Bigr] \Bigr)^{1/q} < \infty.
\]

\section{SDE with rough drift}\label{sec3}

For the statement and proof of the main results we shall use the notion
(and the properties) of
an \textit{SDE with rough drift} captured in the following theorems.
The proofs are postponed to Section~\ref{secRSDEProof}.

As defined above, let $( \bar{\Omega}, \bar{\F}, (\bar{F}_t)_{t\ge
0}, \bar{\P} )$
be a filtered probability space carrying a $d_B$-dimensional Brownian
motion $\bar{B}$ and a bounded $d_S$-dimensional random vector $S_0$
independent of ${\bar{B}}$.
In the following, we fix $\epsilon\in(0,1)$ and $\alpha\in(\frac
{1}{2 + \epsilon}, \frac{1}{2})$. Let $\rp^n\dvtx [0,t] \to\R^{d_Y}$
be smooth paths, such that $\rp^n \to\rrp$ in $\alpha$-H\"older,
for some $\rrp\in\mathcal{C}^{0,\alpha}$, and let $S^n$ be a
$d_S$-dimensional process
which is the unique solution to the classical SDE
\[
S^n_t = S_0 + \int_0^t
a\bigl(S^n_r\bigr) \,dr + \int_0^t
b\bigl(S^n_r\bigr) \,d\bar{B}_r + \int
_0^t c\bigl(S^n_r
\bigr) \,d\rp^n_r,
\]
where we assume that\footnote{In the forthcoming publication
\cite{bibRoughSDE} we show existence of solutions to SDEs with rough
drift under additional sets of assumptions.
The corresponding proofs do not rely on the technique of flow
decomposition used in the present work,
but require more elements of rough path theory and would lead us too
far astray from the topic of filtering.}
%
%
\begin{longlist}[(a1$'$)]
\item[(a1)]
$a \in\Lip^{1}(\R^{d_S})$,
$b_1,\ldots, b_{d_B} \in\Lip^{1}(\R^{d_S})$ and
$c_1,\ldots, c_{d_Y} \in\Lip^{4+\epsilon}(\R^{d_S})$;

\item[(a1$'$)]
$a \in\Lip^{1}(\R^{d_S})$,
$b_1,\ldots, b_{d_B} \in\Lip^{1}(\R^{d_S})$ and
$c_1,\ldots, c_{d_Y} \in\Lip^{5+\epsilon}(\R^{d_S})$.
\end{longlist}

\begin{theorem}
\label{thmRoughSDE} Under assumption \textup{(a1)}, there exists a
$d_S$-dimensional process $S^\infty\in\mathcal{S}^0$ such that
\[
S^n \to S^\infty \qquad\mbox{in } \mathcal{S}^0.
\]
In addition, the limit $\Xi(\rrp):= S^\infty$ only depends on $\rrp
$ and not on the approximating sequence.

Moreover, for all $q\ge1$, $\rrp\in\mathcal{C}^{0,\alpha}$ it
holds that $\Xi(\rrp)
\in\mS^q$ and the corresponding mapping $\Xi\dvtx \mathcal{C}^{0,\alpha
}\to\mathcal
{S}^q$ is locally uniformly continuous [and locally Lipschitz under
assumption \textup{(a1$'$)}].
\end{theorem}

Following Theorem~\ref{thmRoughSDE}, we say that $\Xi(\rrp)$ is a
solution of the \textit{SDE with rough drift}
%
\begin{equation}
\label{eqRough}\qquad \Xi(\rrp)_t = S_0 + \int
_0^t a\bigl(\Xi(\rrp)_r\bigr) \,dr +
\int_0^t b\bigl(\Xi(\rrp)_r
\bigr) \,d\bar{B}_r + \int_0^t c
\bigl(\Xi(\rrp)_r\bigr) \,d\rrp_r.
\end{equation}
The following result establishes some of the salient properties of
solutions of SDEs with rough drift. Recall that $(\Omega, \F, \P_0)$
carries, as above, the $d_Y$-dimensional Brownian motion $Y$,
and let
$\hat{\Omega} = \Omega\times\bar\Omega$ be the product space,
with product measure $\hat{\P}:= \P_0 \otimes\bar{\P}$. Let $S$
be the unique solution on this probability space to the SDE
%
\begin{equation}
\label{eqStratonovich} S_t = S_0 + \int
_0^t a(S_r) \,dr + \int
_0^t b(S_r) \,d\bar{B}_r
+ \int_0^t c(S_r) \circ
dY_r.
\end{equation}
Denote by $\mathbf{Y}$ the rough path lift of $Y$ (i.e., the
enhanced Brownian Motion over $Y$).
%
\begin{theorem}
\label{thmRoughSDE'}
Under assumption \textup{(a1)} we have that:
\begin{itemize}
\item
For every $R > 0, q \ge1$
%
\begin{equation}
\label{eqBoundedExpMoment} \sup_{\|\rrp\|_{\alpha\mbox{-}\Holll} < R} \E\bigl[
\exp\bigl(q \bigl|\Xi(\rrp
)\bigr|_{\infty;[0,t]} \bigr) \bigr] < \infty.
\end{equation}

\item
For $\P_0$-a.e. $\omega$
%
\begin{equation}
\label{eqTheSame} \bar\P\bigl[ S_s(\omega,\cdot) = \Xi\bigl(
\mathbf{Y}(\omega) \bigr)_s(\cdot), s\le t \bigr] = 1.
\end{equation}
\end{itemize}
\end{theorem}

\section{Assumptions and main results}\label{sec4}

In the following we will make use of the Stratonovich version of equation
(\ref{eqMainSDEIto''}); that is, we will consider that the signal
satisfies the equation
%
\begin{eqnarray}
\label{eqMainSDEIto} X_t &=& X_0 + \int
_0^t L_0(X_r,Y_r)
\,dr + \sum_k \int_0^t
Z_k(X_r,Y_r) \circ dY^k_r\nonumber\\
&&{}+ \sum_j \int_0^t
L_j(X_r,Y_r) \,dB^j_r,
\\
Y_t &=& \int_0^t
h(X_r, Y_r) \,dr + W_t,
\nonumber
\end{eqnarray}
where
$L^j_0(x,y) = \bar{l}^j_0(x,y)
- \frac{1}{2} \sum_k \sum_i \partial_{x_i} Z_k^j(x,y) Z_k^i(x,y)
- \frac{1}{2} \sum_k \partial_{y_k} Z_k^j(x,y)$.
We remind the reader that under $\P_0$ the observation $Y$ is a
Brownian motion independent of $B$.

We will assume that $f$ is a bounded Lipschitz function, and we fix
$\epsilon\in(0,1)$
$\alpha\in(\frac{1}{2 + \epsilon}, \frac{1}{2})$,
$t > 0$,
and $X_0$ is a bounded random vector independent of $B$ and $Y$.
We will use one of the following assumptions:
\begin{longlist}[(A1$'$)]
\item[(A1)]
$Z_1,\ldots, Z_{d_Y} \in\Lip^{4 + \epsilon}$,\vspace*{1pt}
$h^1,\ldots, h^{d_Y} \in\Lip^{4 + \epsilon}$ and
$L_0, L_1,\ldots, L_{d_B} \in\Lip^1$;

\item[(A1$'$)]
$Z_1,\ldots, Z_{d_Y} \in\Lip^{5 + \epsilon}$,
$h^1,\ldots, h^{d_Y} \in\Lip^{5 + \epsilon}$ and
$L_0, L_1,\ldots, L_{d_B} \in\Lip^1$.
\end{longlist}
%
\begin{remark}
Assumption (A1) and (A1$'$) lead to the existence of a solution of an
\textit{SDEs with rough driver} (Theorem~\ref{thmRoughSDE}).
Under (A1) the solution mapping is locally uniformly continuous, and
under (A1$'$) it is locally Lipschitz (Theorem~\ref{thmRoughSDE'}).
\end{remark}

Assume either (A1) or (A1$'$).
For $\rrp\in\mathcal{C}^{0,\alpha}$
there exists by Theorem~\ref{thmRoughSDE} a solution $(X^\rrp,
I^{\rrp})$ to the following \textit{SDE with rough drift}:
%
\begin{eqnarray}
\label{eqRoughFilterSDE} X^\rrp_t &=& X_0 +
\int_0^t L_0\bigl(X^\rrp_r,
Y^\rrp_r\bigr) \,dr + \int_0^t
Z\bigl( X^\rrp_r, Y^\rrp_r \bigr) \,d
\rrp_r\nonumber\\
&&{} + \sum_j \int
_0^t L_j\bigl(X^\rrp_r,
Y^\rrp_r\bigr) \,d\bar{B}^j_r,
\nonumber\\[-8pt]\\[-8pt]
Y^\rrp_t &=& \int_0^t d
\rrp_r,
\nonumber\\
I^\rrp_t &=& \int_0^t h
\bigl(X^\rrp_r, Y^\rrp_r\bigr) \,d
\rrp_r - \frac{1}{2} \sum_k \int
_0^t D_k h^k
\bigl(X^\rrp_r, Y^\rrp_r\bigr) \,dr.
\nonumber
\end{eqnarray}

\begin{remark}
Note that formally (!) when replacing the rough path
$\rrp$ with the process $Y$,
$X^\rrp, Y^\rrp$ yields the solution to the SDE (\ref{eqMainSDEIto})
and $\exp( I^\rrp_t )$ yields
the (Girsanov) multiplicator in (\ref{eqPtf''}).
This observation is made precise in the statement of Theorem~\ref{thmRoughSDE}.
\end{remark}

We introduce the functions $g^f,g^1,\theta\dvtx\mathcal{C}^{0,\alpha
}\to
\R$ defined as
\begin{eqnarray*}
g^f(\rrp)&:=& \bar{\E}\bigl[ f\bigl(X^{\rrp}_t,
Y^\rrp_t\bigr) \exp\bigl( I^{\rrp}_t
\bigr) \bigr],\qquad g^1(\rrp):= \bar{\E}\bigl[ \exp\bigl(
I^{\rrp}_t \bigr) \bigr],\\
\theta(\rrp)&:=& \frac{ g^f(\rrp) }{ g^1(\rrp) },\qquad
\rrp\in\mathcal{C}^{0,\alpha}.
\end{eqnarray*}

\begin{theorem}
\label{thmContinuous}
Assume that \textup{(A1)} holds;
then $\theta$ is locally uniformly continuous.
Moreover if \textup{(A1$'$)} holds, then $\theta$ is locally Lipschitz.
\end{theorem}
\begin{pf}
From Theorem~\ref{thmRoughSDE}
we know that for $\rrp\in\mathcal{C}^{0,\alpha}$ the SDE with
rough drift (\ref{eqRoughFilterSDE}) has a unique solution $(X^\rrp,
Y^\rrp, I^\rrp)$
belonging to $\mS^2$.\vadjust{\goodbreak}

Let now $\rrp, \rrp' \in\mathcal{C}^{0,\alpha}$.
Denote $X = X^{\rrp}, Y = Y^{\rrp}, I = I^{\rrp}$ and analogously
for~$\rrp'$.

Then
\begin{eqnarray*}
&&
\bigl|g^f(\rrp) - g^f(\rrp)\bigr| \\
&&\qquad\le\E\bigl[ \bigl|f(X_t,Y_t) \exp( I_t ) - f \bigl(X'_t,Y"_t\bigr)
\exp\bigl( I'_t \bigr)\bigr| \bigr]
\\
&&\qquad\le\E\bigl[ \bigl|f(X_t,Y_t)\bigr| \bigl|\exp( I_t ) - \exp\bigl( I'_t \bigr)\bigr|
\bigr] \\
&&\qquad\quad{} + \E \bigl[ \bigl|f(X_t,Y_t) - f\bigl(X'_t,Y'_t \bigr)\bigr| \exp\bigl(
I'_t \bigr) \bigr]
\\
&&\qquad\le|f|_{\infty} \E\bigl[ \bigl|\exp( I_t ) - \exp \bigl( I'_t \bigr)\bigr|
\bigr]\\
&&\qquad\quad{} + \E\bigl[ \bigl|f(X_t,Y_t) - f\bigl(X'_t,Y'_t \bigr)\bigr|^2 \bigr]^{1/2}
\E\bigl[ \bigl|\exp\bigl( I'_t \bigr)\bigr|^2 \bigr]^{1/2}
\\
&&\qquad\le |f|_{\infty} \E\bigl[ \bigl|\exp( I_t ) + \exp
\bigl( I'_t \bigr)\bigr|^2 \bigr]^{1/2}
\E\bigl[ \bigl|I_t - I'_t\bigr| \bigr]^{1/2}
\\
&&\qquad\quad{}+ \E\bigl[ \bigl|f(X_t,Y_t) - f\bigl(X'_t,Y'_t
\bigr)\bigr|^2 \bigr]^{1/2} \E\bigl[ \bigl|\exp\bigl(
I'_t \bigr)\bigr|^2 \bigr]^{1/2}.
\end{eqnarray*}

Hence, using from Theorems~\ref{thmRoughSDE} and~\ref{thmRoughSDE'} the
continuity statements
as well as the boundedness of exponential moments,
we see that $g^f$ is locally uniformly continuous under (A1),
and it is locally Lipschitz under (A1$'$).

The same then holds true for $g^1$ and
moreover $g^1(\rrp) > 0$.
Hence $\theta$ is locally uniformly continuous under (A1)
and
locally Lipschitz under (A1$'$).
\end{pf}

Denote by $\mathbf{Y}_\cdot$, as before, the canonical rough path
lift of $Y$ to $\mathcal{C}^{0,\alpha}$.
We then have
%
\begin{theorem}
\label{thmVersionOfConditionalExpectation}
Assume either \textup{(A1)} or \textup{(A1$'$)}. Then
$ \theta( \mathbf{Y}_\cdot) = \pi_t(f)$,
$\P$-a.s.
\end{theorem}
\begin{pf}
To prove the statement it is enough to show that
\[
g^f( \mathbf{Y}_\cdot) = p_t(f),\qquad \P
\mbox{-a.s.},
\]
which is equivalent to
\[
g^f( \mathbf{Y}_\cdot) = p_t(f),\qquad
\P_0\mbox{-a.s.}
\]

For that, it suffices to show that
%
\begin{equation}
\label{eqEqualityToShow} \E_0\bigl[ p_t(f) \Upsilon(
Y_\cdot) \bigr] = \E_0\bigl[ g^f( \mathbf
{Y}_\cdot) \Upsilon( Y_\cdot) \bigr]
\end{equation}
for an arbitrary continuous bounded function $\Upsilon\dvtx C([0,t], \R
^{d_Y}) \to\R$.

Let $( \bar{\Omega}, \bar{\F}, \bar{\P} )$ be the auxiliary
probability space from before,
carrying an $d_B$-dimensional Brownian motion $\bar{B}$.
Let $(\hat{\Omega}, \hat{\F}, \hat{\P}):= (\Omega\times\bar
{\Omega}, \F\otimes\bar{\F}, \P_0 \otimes\bar{\P})$.
By $Y$ and $X_0$ we denote also the ``lift'' of $Y$ to $\hat{\Omega}$,
that is, $Y(\omega,\bar{\omega}) = Y(\omega)$, $X_0(\omega,\bar
{\omega}) = X_0(\omega)$.
Then $(Y, B)$ (on $\Omega$ under $\P_0$) has the same distribution as
$(Y, \bar{B})$ (on $\hat{\Omega}$ under $\hat{\P}$).

Denote by $(\hat{X}, \hat{I})$ the solution on $(\hat{\Omega}, \hat
{F}, \hat{\P})$ to the SDE
\begin{eqnarray*}
\hat{X}_t &=& X_0 + \int_0^t
L_0(\hat{X}_r, Y_r) \,dr + \sum
_k \int_0^t
Z_k( \hat{X}_r, Y_r) \circ
dY^k_r \\
&&{}+ \sum_j \int
_0^t L_j(\hat{X}_r,
Y_r) \,d\bar{B}^j_r,
\\
\hat{I}_t &=& \sum_k \int
_0^t h^k(\hat{X}_r,
Y_r) \circ dY^k_r - \frac{1}{2} \sum
_k \int_0^t
D_k h^k(\hat{X}_r, Y_r) \,dr.
\end{eqnarray*}

Then
\[
(Y, \hat{X}, \hat{I})_{\hat{\P}} \sim\biggl(Y, X, \sum
_k \int_0^\cdot
h^k(X_r, Y_r) \circ dY^k_r
- \frac{1}{2} \sum_k \int
_0^\cdot D_k h^k(X_r,
Y_r) \,dr \biggr)_{\P_0}.
\]

Hence, for the left-hand side of (\ref{eqEqualityToShow}),
\begin{eqnarray*}
&&
\E_0\bigl[ p_t(f) \Upsilon( Y_\cdot) \bigr] \\
&&\qquad=
\E_0\biggl[ f( X_t, Y_t ) \exp\biggl( \sum
_k \int_0^t
h^k(X_r,Y_r) \circ dY^k_r\\
&&\qquad\quad\hspace*{82.2pt}{} - \frac{1}{2} \sum_k \int
_0^t D_k h^k(X_r,Y_r)
\,dr \biggr) \Upsilon( Y_\cdot) \biggr]
\\
&&\qquad= \hat{\E}\bigl[ f( \hat{X}_t, Y_t ) \exp(
\hat{I}_t ) \Upsilon(Y_\cdot) \bigr].
\end{eqnarray*}

On the other hand, from Theorem~\ref{thmRoughSDE'} we know that
for $\P_0$-a.e. $\omega$
\begin{eqnarray*}
X^{ \mathbf{Y}_\cdot(\omega) }(\bar{\omega})_t &=& \hat{X}_t(\omega,
\bar{\omega}),\qquad Y^{ \mathbf{Y}_\cdot(\omega) }(\bar{\omega})_t = \hat
{Y}_t(\omega, \bar{\omega}),\\
I^{ \mathbf{Y}_\cdot(\omega) }(\bar{
\omega})_t &=& \hat{I}_t(\omega, \bar{\omega}), \qquad\bar{\P}
\mbox{-a.e. } \bar{\omega}. 
\end{eqnarray*}

Hence, for the right-hand side of (\ref{eqEqualityToShow}) we get
(using Fubini for the last equality)
\begin{eqnarray*}
\E_0\bigl[ g^f( \mathbf{Y}_\cdot)
\Upsilon( Y_\cdot) \bigr] &=& \E_0\bigl[ \bar{\E}\bigl[ f
\bigl( X^{\mathbf{Y}_\cdot}_t, Y^{\mathbf
{Y}_\cdot}_t \bigr) \exp
\bigl( I^{\mathbf{Y}_\cdot}_t \bigr) \bigr] \Upsilon( Y_\cdot)
\bigr]
\\
&=& \E_0\bigl[ \bar{\E}\bigl[ f( \hat{X}_t,
Y_t ) \exp( \hat{I}_t ) \bigr] \Upsilon(
Y_\cdot) \bigr]
\\
&=& \hat{\E}\bigl[ f( \hat{X}_t, Y_t ) \exp(
\hat{I}_t ) \Upsilon( Y_\cdot) \bigr],
\end{eqnarray*}
which yields (\ref{eqEqualityToShow}).
\end{pf}

\section{\texorpdfstring{Proofs of Theorems \protect\ref{thmRoughSDE} and \protect\ref{thmRoughSDE'}}
{Proofs of Theorems 2 and 3}}\label{sec5} \label{secRSDEProof}

\mbox{}

\begin{pf*}{Proof of Theorem~\ref{thmRoughSDE}}
Let $\rrp\in\mathcal{C}^{0,\alpha}$ be the lift of a smooth path
$\rp$.
Let $S^\rrp$ be the unique solution of the SDE
\[
S^\rrp_t = S_0 + \int_0^t
a\bigl(S^\rrp_r\bigr) \,dr + \int_0^t
b\bigl(S^\rrp_r\bigr) \,d\bar{B}_r + \int
_0^t c\bigl(S^\rrp_r
\bigr) \,d\rp_r.
\]

Define $\tilde{S}^\rrp:= (\olphi^\rrp)^{-1}(t,S^\rrp_t)$, where
$\olphi^\rrp$ is the ODE flow
%
\begin{equation}
\label{eqODE} \olphi^\rrp(t,x) = x + \int_0^t
c\bigl( \olphi^\rrp(r,x) \bigr) \,d\rrp_r.
\end{equation}

By Lemma~\ref{lemTransformation}, we have that $\tilde{S}^\rrp$
satisfies the SDE
%
\begin{equation}
\label{eqSDETransformed} \tilde{S}^\rrp_t =
S_0 + \int_0^t
\tilde{a}^\rrp\bigl(r, \tilde{S}^\rrp_r\bigr) \,dr
+ \int_0^t \tilde{b}^\rrp\bigl(r,
\tilde{S}^\rrp_r\bigr) \,d\bar{B}_r
\end{equation}
with $\tilde{a}^\rrp, \tilde{b}^\rrp$ defined as in Lemma \ref
{lemTransformation}.

This equation makes sense, even if $\rrp$ is a generic rough path in
$\mathcal{C}^{0,\alpha}$
[in which case (\ref{eqODE}) is now really an RDE].
Indeed, since the first two derivatives of $\olphi^\rrp$ and its
inverse are bounded
(Proposition 11.11 in~\cite{bibFrizVictoir})
we have that $\tilde{a}^\rrp(t,\cdot), \tilde{b}^\rrp(t,\cdot)$
are also in $\Lip^1$.
Hence by Theorem V.7 in~\cite{bibProtter},
there exists a unique strong solution to (\ref{eqSDETransformed}).

We define the mapping introduced in Theorem~\ref{thmRoughSDE} as
\[
\Xi(\rrp)_t:= \olphi^\rrp\bigl(t, \tilde{S}^\rrp_t
\bigr).
\]

To show continuity of the mapping we restrict ourselves to the case $q=2$.
Moreover we shall assume
$c_1,\ldots, c_{d_Y} \in\Lip^{5+\epsilon}(\R^{d_S})$, and we will
hence prove the local Lipschitz property of the respective maps.

Let $\rrp^1, \rrp^2 \in\mathcal{C}^{0,\alpha}$ with $|\rrp^1|_{\alpha
\mbox{-}\Holl}, |\rrp^2|_{\alpha\mbox{-}\Holl} < R$.
By Lemma~\ref{lemSDELipschitzDependence} we have
\[
\bar{\E}\Bigl[ \sup_{s\le t} \bigl|\tilde{S}^1_s -
\tilde{S}^2_s\bigr|^2 \Bigr]^{1/2} \le
C_{\mmbox{Lem~\ref{lemSDELipschitzDependence}}}(R) \rho_{\alpha\mbox{-}\Holl}\bigl(\rrp^1,
\rrp^2\bigr).
\]

Hence
\begin{eqnarray*}
&&
\bar{\E}\Bigl[ \sup_{s\le t} \bigl|\Xi\bigl(\rrp^1
\bigr)_s - \Xi\bigl(\rrp^2\bigr)_s\bigr|^2
\Bigr]^{1/2} \\
&&\qquad= \bar{\E}\Bigl[ \sup_{s\le t} \bigl|
\phi^1\bigl(s,\tilde{S}^1_s\bigr) -
\phi^2\bigl(s,\tilde{S}^2_s
\bigr)\bigr|^2 \Bigr]^{1/2}
\\
&&\qquad\le\bar{\E}\Bigl[ \sup_{s\le t} \bigl|\phi^1\bigl(s,
\tilde{S}^1_s\bigr) - \phi^1\bigl(s,
\tilde{S}^2_s\bigr)\bigr|^2 \Bigr]^{1/2}
+ \bar{\E}\Bigl[ \sup_{s\le t} \bigl|\phi^1\bigl(s,
\tilde{S}^2_s\bigr) - \phi^2\bigl(s,
\tilde{S}^2_s\bigr)\bigr|^2 \Bigr]^{1/2}
\\
&&\qquad\le\bar{\E}\Bigl[ \sup_{s\le t} \bigl|\phi^1\bigl(s,
\tilde{S}^1_s\bigr) - \phi^1\bigl(s,
\tilde{S}^2_s\bigr)\bigr|^2 \Bigr]^{1/2}
+ \sup_{s\le t, x \in\R^{d_S} } \bigl|\phi^1(s,x) - \phi^1(s,x)\bigr|
\\
&&\qquad\le K(R) \bar{\E}\Bigl[ \sup_{s\le t} \bigl|\tilde{S}^1_s
- \tilde{S}^2_s\bigr|^2 \Bigr]^{1/2} +
C_{\mmbox{Lem~\ref{lemLipschitnessOfFlow}}}(R) \rho_{\alpha\mbox{-}\Holl}\bigl(\rrp^1,
\rrp^2\bigr)
\\
&&\qquad\le C_1 \rho_{\alpha\mbox{-}\Holl}\bigl(\rrp^1,
\rrp^2\bigr)
\end{eqnarray*}
as desired, where
\[
C_1 = K_{\mmbox{Lem~\ref{lemLipschitnessOfFlow}}}(R) C_{\mmbox{Lem \ref
{lemSDELipschitzDependence}}}(R) + C_{\mmbox{Lem~\ref{lemLipschitnessOfFlow}}}(R),
\]
where $K_{\mmbox{Lem~\ref{lemLipschitnessOfFlow}}}(R)$ and $C_{\mmbox{Lem \ref
{lemLipschitnessOfFlow}}}(R)$ are the constants from Lemma \ref
{lemLipschitnessOfFlow}.\vadjust{\goodbreak}
\end{pf*}
\begin{pf*}{Proof of Theorem~\ref{thmRoughSDE'}}
In order to show (\ref{eqBoundedExpMoment}),
pick $k \in1,\ldots, d_S$.
We first note that by simply scaling (the coefficients of) $S^\rrp$ it
is sufficient to argue for $q=1$.
And consider the $k$th component of $\Xi$.

Then
\begin{eqnarray*}
&&
\E\bigl[ \exp\bigl( \bigl|\Xi^{(k)}(\rrp)\bigr|_{\infty;[0,t]} \bigr) \bigr]
\\
&&\qquad\le\E\bigl[ \exp\bigl( \bigl|D \olpsi^{\rrp}\bigr|_\infty
\bigl( \bigl|\olphi^{\rrp}(0,S_0)\bigr| + \bigl|\tS^{(k);\rrp}\bigr|_{\infty;[0,t]}
\bigr) \bigr) \bigr]
\\
&&\qquad\le\exp\Bigl( \bigl|D \olpsi^{\rrp}\bigr|_\infty
\sup_{|x|\le|S_0|_{L^\infty}} \bigl|\olphi^{\rrp}(0, x)\bigr| \Bigr) \E\Bigl[ \exp
\Bigl( \bigl|D
\olpsi^{\rrp}\bigr|_\infty\sup_{s\le t} \bigl|\tS^{(k);\rrp
}_t\bigr|
\Bigr) \Bigr]
\\
&&\qquad\le\exp\Bigl( \bigl|D \olpsi^{\rrp}\bigr|_\infty
\sup_{|x|\le|S_0|_{L^\infty}} \bigl|\olphi^{\rrp}(0, x)\bigr| \Bigr)
\\
&&\qquad\quad{} \times\Bigl( \E\Bigl[ \exp\Bigl( \bigl|D \olpsi^{\rrp}\bigr|_\infty
\sup_{s\le t} \tS^{(k);\rrp
}_s\Bigr) \Bigr] + \E\Bigl[ -
\exp\Bigl( \bigl|D \olpsi^{\rrp}\bigr|_\infty\sup_{s\le t}
\tS^{(k);\rrp
}_s\Bigr) \Bigr] \Bigr)
\\
&&\qquad= \exp\Bigl( \bigl|D \olpsi^{\rrp}\bigr|_\infty
\sup_{|x|\le|S_0|_{L^\infty}} \bigl|\olphi^{\rrp}(0, x)\bigr| \Bigr)
\\
&&\qquad\quad{} \times\Bigl( \E\Bigl[ \sup_{s\le t} \exp\bigl( \bigl|D \olpsi^{\rrp
}\bigr|_\infty
\tS^{(k);\rrp
}_s\bigr) \Bigr] + \E\Bigl[ - \sup_{s\le t}
\exp\bigl( \bigl|D \olpsi^{\rrp}\bigr|_\infty\tS^{(k);\rrp
}_s
\bigr) \Bigr] \Bigr).
\end{eqnarray*}
Now, only the boundedness of the last two terms remains to be shown,
for $\rrp$ bounded.

By applying It\^{o}'s formula we get that
\begin{eqnarray*}
\exp\bigl( \tS^{(k);\rrp}_t \bigr) &=& 1 + \int
_0^t \exp\bigl( \tS^{(k);\rrp}_r
\bigr) \,d\tS^{(k);\rrp}_r + \int_0^t
\exp\bigl( \tS^{(k);\rrp}_r \bigr) \,d\bigl\langle
\tS^{(k);\rrp} \bigr\rangle_r
\\
&=& 1 + \int_0^t \exp\bigl(
\tS^{(k);\rrp}_r \bigr) \tilde{a}^{\rrp}_k
\bigl( \tS^{(k);\rrp}_r \bigr) \,dr + \sum
_{i=1}^{d_B} \int_0^t
\exp\bigl( \tS^{(k);\rrp}_r \bigr) \tilde{b}^{\rrp}_{k i}
\bigl( \tS^{(k);\rrp}_r \bigr) \,d\bar{B}^i_r
\\
&&{} + \sum_{i=1}^{d_B} \int
_0^t \exp\bigl( \tS^{(k);\rrp}_r
\bigr) \bigl|\tilde{b}^{\rrp}_{k i}\bigl( \tS^{(k);\rrp}_r
\bigr)\bigr|^2 \,dr.
\end{eqnarray*}
Hence the process $\exp( \tS^{(k);\rrp} )$ satisfies an SDE with
Lipschitz coefficients and by
an application of Gronwalls lemma and the Burkholder--Davis--Gundy inequality
(see also Lemma V.2 in~\cite{bibProtter})
one arrives at
%
\begin{equation}
\label{eqExpFinite} \sup_{|\rrp|_{\alpha\mbox{-}\Holll} < R} \sup_{s\le t} \bar
{\E}\bigl[ \exp
\bigl( \bigl|D \olpsi^{\rrp}\bigr|_\infty\tS^{(k);\rrp}_t
\bigr) \bigr] \le C \exp( C_2 ),
\end{equation}
where $C$ is universal and
\[
C_2:= \sup_{ \rrp: \|\rrp\|_{\alpha\mbox{-}\Holll} < R } \bigl|\tilde{a}^\rrp
\bigr|_\infty+ \bigl|\tilde{b}^\rrp\bigr|_\infty^2,
\]
which is finite because of Lemma~\ref{lemTransformationRough}.
One argues analogously for
\[
\sup_{s\le t} \bar{\E}\bigl[ - \exp\bigl( \bigl|D \olpsi^{\rrp}\bigr|_\infty\tS^{(k);\rrp
}_t \bigr) \bigr],
\]
which then gives (\ref{eqBoundedExpMoment}).

Now, for the correspondence to an SDE solution let $\Omega$ be the additional
probability space as given in the statement.
Let $S$ be the solution to the SDE
(\ref{eqStratonovich}).

In Section 3 in~\cite{bibBuckdahnMaHJB}
it was shown (see also Theorem 2 in~\cite{bibBenArousCastell}), that
if we let $\Theta$ be the stochastic (Stratonovich) flow
\[
\Theta(\omega;t,x) = x + \int_0^t c\bigl(
\Theta(\omega;r,x) \bigr) \circ dY_r(\omega),
\]
then with $\hat{S}_t:= \Theta^{-1}(t, S_t)$ we have $\hat{\P}$-a.s.
%
\begin{eqnarray}
\label{eqSDETransformedStochastic} \hat{S}_s\bigl(\omega,
\omega^{\bar{B}}\bigr) &=& S_0 + \int_0^s
\hat{a}(r, \hat{S}_r) \,dr\nonumber\\[-9pt]\\[-9pt]
&&{} + \int_0^s
\hat{b}(r, \hat{S}_r) \,d\bar{B}_r,\qquad s \in[0,t], \qquad\hat{\P}
\mbox{-a.e. } \bigl(\omega,\omega^{\bar B}\bigr).\nonumber
\end{eqnarray}
Here, componentwise,
\begin{eqnarray*}
\hat{a}(t,x)_i &:=& \sum_k
\partial_{x_k} \Theta^{-1}_i\bigl(t,\Theta(t,x)
\bigr) a_k\bigl( \Theta(t,x) \bigr) \\[-2pt]
&&{}+ \frac{1}{2} \sum
_{j,k} \partial_{x_j x_k} \Theta^{-1}_i
\bigl(t,\Theta(t,x)\bigr) \sum_l
b_{jl}\bigl(\Theta(t,x)\bigr) b_{kl}\bigl(\Theta(t,x)
\bigr),
\\[-2pt]
\hat{b}(t,x)_{ij} &:=& \sum_k
\partial_{x_k} \Theta^{-1}_i(t,x)
b_{k
j}\bigl(\Theta(t,x)\bigr).
\end{eqnarray*}

Especially, by a Fubini-type theorem (e.g., Theorem 3.4.1 in \cite
{bibBogachev}), there exists $\Omega_0$ with $\P_0(\Omega_0) = 1$
such that for $\omega\in\Omega_0$ equation (\ref
{eqSDETransformedStochastic}) holds true $\bar\P$-a.s.

Let $\mathbf{Y} \in\mathcal{C}^{0,\alpha}$ be the enhanced Brownian
motion over $Y$.
We can then construct $\omega$-wise the rough flow $\phi^{ \mathbf
{Y}(\omega) }$
as given in (\ref{eqODE}).
By the very definition of $\Xi$ we know that
$\tilde{S}^{\mathbf{Y}(\omega)}_t(\omega):= (\phi^{\mathbf
{Y}(\omega)})^{-1}(\omega;t,\Xi(\bolds{\omega})_t)$
satifies the SDE
%
\begin{eqnarray}
\label{eqSDETildeS} \tilde{S}^{\mathbf{Y}(\omega)}_t &=& S_0 +
\int_0^t \hat{b}^{\mathbf
{Y}(\omega)}\bigl(r,
\tilde{S}^{\mathbf{Y}(\omega)}_r\bigr) \,dr\nonumber\\[-9pt]\\[-9pt]
&&{} + \int_0^t
\hat{b}^{\mathbf{Y}(\omega)}\bigl(r,\tilde{S}^{\mathbf{Y}(\omega
)}_r\bigr) \,d
\bar{B}_r, \qquad\bar\P\mbox{-a.e. }\omega^{\bar B},\nonumber
\end{eqnarray}
where
\begin{eqnarray*}
\tilde{a}^{\mathbf{Y}(\omega)}(t,x)_i &:=& \sum
_k \partial_{x_k} \bigl(\phi^{\mathbf{Y}(\omega)}
\bigr)^{-1}_i\bigl(t,\phi^{\mathbf{Y}(\omega)}(t,x)\bigr)
a_k\bigl(t, \phi^{\mathbf{Y}(\omega)}(t,x) \bigr)
\\
&&{} + \frac{1}{2} \sum_{j,k}
\partial_{x_j x_k} \bigl(\phi^{\mathbf{Y}(\omega
)}\bigr)^{-1}_i
\bigl(t,\phi^{\mathbf{Y}(\omega)}(t,x)\bigr) \\
&&\hspace*{35.5pt}{}\times\sum_l
b_{jl}\bigl(t,\phi^{\mathbf{Y}(\omega)}(t,x)\bigr) b_{kl}\bigl(t,
\phi^{\mathbf{Y}(\omega
)}(t,x)\bigr),
\\
\tilde{b}^{\mathbf{Y}(\omega)}(t,x)_{ij} &:=& \sum
_k \partial_{x_k} \bigl(\phi^{\mathbf{Y}(\omega)}
\bigr)^{-1}_i(t,x) b_{k
j}\bigl(t,
\phi^{\mathbf{Y}(\omega)}(t,x)\bigr).
\end{eqnarray*}

It is a classical rough path result
(see, e.g., Section 17.5 in~\cite{bibFrizVictoir}),
that there exists $\Omega_1$ with $\P^Y(\Omega_1) = 1$ such that for
$\omega\in\Omega_1$, we have
\[
\phi^{\mathbf{Y}(\omega)}(\cdot,\cdot) = \Theta(\omega;\cdot,\cdot).
\]
Hence for $\omega\in\Omega_1$ we have that $\hat{a} = \tilde
{a}^{\mathbf{Y}(\omega)}, \hat{b} = \tilde{b}^{\mathbf{Y}(\omega)}$.
Hence for \mbox{$\omega\in\Omega_0 \cap\Omega_1$} the processes
$\hat{S}_t(\omega,\cdot), \tilde{S}^{\mathbf{Y}(\omega)}_t(\cdot)$
satisfy the same Lipschitz SDE (with respect
to~$\bar\P$).\footnote{Here one has to argue
that fixing $\omega$ in equation (\ref{eqSDETildeS})
gives ($\P_0$-a.s.) the solution to the respective SDE on $\Omega^{\bar B}$.}
By strong uniqueness we hence have
for $\omega\in\Omega_0 \cap\Omega_1$ that $\bar\P$-a.s.
\[
\hat{S}_s(\omega,\cdot) = \tilde{S}^{\mathbf{Y}(\omega)}_s(
\cdot),\qquad s\le t.
\]

Hence for $\omega\in\Omega_0 \cap\Omega_1$
\[
S_s(\omega,\cdot) = \Xi\bigl(\mathbf{Y}(\omega)\bigr) (
\cdot)_s,\qquad s\le t,\qquad \bar\P\mbox{-a.s.}
\]
\upqed
\end{pf*}
%
\begin{remark}
We remark that the above idea of a flow decomposition
is also used in the work by Davis~\cite{d1,d2,d3,d4}.
Without rough path theory this approach is restricted to
one-dimensional observation,
since, for multidimensional flows, one cannot hope for continuous
dependence on the driving signal
in supremum norm.
\end{remark}
%
\begin{lemma}
\label{lemTransformation}
Let $\rp$ be a smooth $d_Y$-dim path $\rp$ and $S$ be the solution of
the the following classical SDE
\[
S_t = S_0 + \int_0^t
a(S_r) \,dr + \int_0^t
b(S_r) \,d\bar{B}_r + \int_0^t
c(S_r) \,d\rp_r,
\]
where $\bar{B}$ is a $d_B$-dimensional Brownian motion,
\[
\int_0^t c(S_r) \,d\rp_r:= \sum_{i=1}^{d_S} \int_0^t c_i(S_r) \dot
{\rp}^i_r \,dr,
\]
$a \in\Lip^{1}(\R^{d_S})$,
$b_1,\ldots, b_{d_B} \in\Lip^1( \R^{d_S} )$,
$c_1,\ldots, c_{d_Y} \in\Lip^{4 + \epsilon}(\R^{d_S})$,
and $S_0 \in L^\infty(\bar\Omega; \R^{d_S})$ independent of $\bar{B}$.
Consider the flow
%
\begin{equation}
\label{eqOlphi} \olphi(t,x) = x + \int_0^t c
\bigl( \olphi(r,x) \bigr) \,d\rp_r.
\end{equation}

Then $\tilde{S}_t:= \phi^{-1}(t, S_t)$ satisfies the following SDE:
\[
\tilde{S}_t = S_0 + \int_0^t
\tilde{a}(r,\tilde{S}_r) \,dr + \int_0^t
\tilde{b}(r,\tilde{S}_r) \,d\bar{B}_r,
\]
where we define componentwise
\begin{eqnarray*}
\tilde{a}(t,x)_i &:=& \sum_k
\partial_{x_k} \olphi^{-1}_i\bigl(t,\olphi(t,x)
\bigr) a_k\bigl(\olphi(t,x) \bigr)\\
&&{} + \frac{1}{2} \sum
_{j,k} \partial_{x_j x_k} \olphi^{-1}_i
\bigl(t, \olphi(t,x)\bigr) \sum_l
b_{jl}\bigl(\olphi(t,x)\bigr) b_{kl}\bigl(\olphi(t,x)
\bigr),
\\
\tilde{b}(t,x)_{ij} &:=& \sum
_k \partial_{x_k} \olphi^{-1}_i
\bigl(t,\olphi(t,x)\bigr) b_{k j}\bigl(\olphi(t,x)\bigr).
\end{eqnarray*}
\end{lemma}
\begin{pf}
Denote $\olpsi(t,x):= \olphi^{-1}(t,x)$. Then
\[
\olpsi(r,x) = x - \int_0^t
\partial_x \olpsi(r,x) c(x) \,d\rp_r.
\]

By It\^o's formula,
\begin{eqnarray*}
&&
\olpsi_i(t, S_t) - \olpsi_i(0,S_0)\\
&&\qquad= \int_0^t \partial_t
\olpsi_i(r,S_r) \,dr + \sum_j
\int_0^t \partial_{x_j}
\olpsi_i(r,S_r) \,dS_j(r) \\
&&\qquad\quad{}+ \sum
_{j,k} \frac{1}{2} \int_0^t
\partial_{x_j x_k} \olpsi_i(r,S_r) \,d\langle
S_k, S_j \rangle_r
\\
&&\qquad= \sum_j \int_0^t
\partial_{x_j} \olpsi_i(r,S_r)
a_j(S_r) \,dr + \sum_j \int
_0^t \partial_{x_j}
\olpsi_i(r,S_r) \sum_k
b_{jk}(S_r) \,d\bar{B}_k(r)
\\
&&\qquad\quad{} + \sum_{j,k} \frac{1}{2} \int
_0^t \partial_{x_j x_k}
\olpsi_i(S_r) \sum_l
b_{kl}(S_r) b_{jl}(S_r) \,dr.
\end{eqnarray*}
\upqed
\end{pf}
%
\begin{lemma}
\label{lemTransformationRough}
Consider for a rough path $\rrp\in\mathcal{C}^{0,\alpha}$
the coefficients transformed analogously to Lemma \ref
{lemTransformation}, $\tilde{a}^\rrp$, $\tilde{b}^\rrp$;
that is, consider the rough flow
%
\begin{equation}
\label{eqOlphiRough} \olphi(t,x) = \olphi^\rrp(t,x) = x + \int
_0^t c\bigl( \olphi(r,x) \bigr) \,d
\rrp_r
\end{equation}
and define
\begin{eqnarray*}
\tilde{a}^\rrp(t,x)_i &:=& \sum
_k \partial_{x_k} \olphi^{-1}_i
\bigl(t,\olphi(t,x)\bigr) a_k\bigl(\olphi(t,x) \bigr)\\
&&{} +
\frac{1}{2} \sum_{j,k} \partial_{x_j x_k}
\olphi^{-1}_i\bigl(t, \olphi(t,x)\bigr) \sum
_l b_{jl}\bigl(\olphi(t,x)\bigr)
b_{kl}\bigl(\olphi(t,x)\bigr),
\\
\tilde{b}^\rrp(t,x)_{ij} &:=& \sum
_k \partial_{x_k} \olphi^{-1}_i
\bigl(t,\olphi(t,x)\bigr) b_{k j}\bigl(\olphi(t,x)\bigr).
\end{eqnarray*}

Then for every $R > 0$ there exists $K_{\mmbox{Lem \ref
{lemTransformationRough}}} = K_{\mmbox{Lem~\ref{lemTransformationRough}}}(R) <
\infty$ such that
\begin{eqnarray*}
\sup_{\rrp: |\rrp|_{\alpha\mbox{-}\Holll} < R} \bigl|\tilde{a}^\rrp\bigr|_\infty&\le&
K_{\mmbox{Lem~\ref{lemTransformationRough}}},
\\
\sup_{\rrp: |\rrp|_{\alpha\mbox{-}\Holll} < R} \bigl|\tilde{b}^\rrp\bigr|_\infty&\le&
K_{\mmbox{Lem~\ref{lemTransformationRough}}},
\\
\sup_{\rrp: |\rrp|_{\alpha\mbox{-}\Holll} < R} \sup_{s\le t}\bigl|D \tilde{a}^\rrp(s,
\cdot)\bigr|_\infty&\le& K_{\mmbox{Lem~\ref{lemTransformationRough}}},
\\
\sup_{\rrp: |\rrp|_{\alpha\mbox{-}\Holll} < R} \sup_{s\le t}\bigl|D \tilde{b}^\rrp(s,
\cdot)\bigr|_\infty&\le& K_{\mmbox{Lem~\ref{lemTransformationRough}}},
\end{eqnarray*}
and such that if $\rrp, \tilde\rrp$ are two rough paths with $|\rrp
|_{\alpha\mbox{-\textup{\Holl}}}, |\tilde\rrp|_{\alpha\mbox{-\textup{\Holl}}} < R$, we have
\begin{eqnarray*}
&&\sup_{t,x} \bigl|\tilde{a}^1(t,x) - \tilde{a}^2(t,x)\bigr|
\le K_{\mmbox{Lem~\ref{lemTransformationRough}}}(R) \rho_{\alpha\mbox{-\textup{\Holl}}}\bigl(
\rrp^1,
\rrp^2 \bigr),
\\
&&\sup_{t,x} \bigl|\tilde{b}^1(t,x) - \tilde{b}^2(t,x)\bigr|
\le K_{\mmbox{Lem~\ref{lemTransformationRough}}}(R) \rho_{\alpha\mbox{-\textup{\Holl}}}\bigl(
\rrp^1,
\rrp^2 \bigr).
\end{eqnarray*}
\end{lemma}
\begin{pf}
This is a straightforward calculation using
Lemma~\ref{lemLipschitnessOfFlow} and
the properties of $a, b$.
\end{pf}

The following is a standard result for continuous dependence of SDEs on
parameters.
%
\begin{lemma}
\label{lemSDELipschitzDependencePureIto}
Let $\tilde{a}^i(t,x)$, $\tilde{b}^i(t,x), i=1,2$, be bounded and
uniformly Lipschitz in $x$.

Let $\tilde{S}^i$ be the corresponding unique solutions to the SDEs
\[
\tilde{S}^i_t = S_0 + \int
_0^t \tilde{a}^i\bigl( r,
\tilde{S}^i_r \bigr) \,dr + \int_0^t
\tilde{b}^i\bigl( r, \tilde{S}^i_r \bigr) \,d
\bar{B}_r,\qquad i=1,2.
\]

Assume
\begin{eqnarray*}
&\displaystyle \sup_{s\le t} \bigl|D \tilde{a}^1(s,\cdot)\bigr|_\infty,\qquad
\sup_{s\le t} \bigl|D \tilde{b}^1(s,\cdot)\bigr|_\infty< K <
\infty,&
\\
&\displaystyle \sup_{r,x} \bigl(\tilde{a}^1(r,x) - \tilde{a}^2(r,x)
\bigr),\qquad \sup_{r,x} \bigl(\tilde{b}^1(r,x) -
\tilde{b}^2(r,x) \bigr) < \varepsilon< \infty.&
\end{eqnarray*}

Then there exists $C_{\mmbox{Lem~\ref{lemSDELipschitzDependencePureIto}}} =
C_{\mmbox{Lem~\ref{lemSDELipschitzDependencePureIto}}}( K )$
such that
\[
\E\Bigl[ \sup_{ s\le t} \bigl|\tilde{S}^1_s -
\tilde{S}^2_s\bigr|^2 \Bigr]^{1/2} \le
C_{\mmbox{Lem~\ref{lemSDELipschitzDependencePureIto}}} \varepsilon.
\]
\end{lemma}
\begin{pf}
This is a straightforward application of It\^{o}'s formula and
the
Burkholder--Davis--Gundy inequality.
\end{pf}

We now apply the previous lemma to our concrete setting.
%
\begin{lemma}
\label{lemSDELipschitzDependence}
Let $\rrp^1, \rrp^2 \in\mathcal{C}^{0,\alpha}$
and let $\tilde{S}^1, \tilde{S}^2$ be the corresponding unique solutions
to the SDEs
\[
\tilde{S}^i_t = S_0 + \int
_0^t \tilde{a}^i\bigl( r,
\tilde{S}^i_r \bigr) \,dr + \int_0^t
\tilde{b}^i\bigl( r, \tilde{S}^i_r \bigr) \,d
\bar{B}_r,\qquad i=1,2,
\]
where $\tilde{a}^i, \tilde{b}^i$ are given as in Lemma \ref
{lemTransformationRough}.

Assume $R > \max\{ |\rrp^1|_{\alpha\mbox{-\textup{\Holl}}},
|\rrp^2|_{\alpha\mbox{-\textup{\Holl}}
} \}$.
Then there exists $C_{\mmbox{Lem~\ref{lemSDELipschitzDependence}}} =\break C_{\mmbox{Lem
\ref{lemSDELipschitzDependence}}}( R )$:
\[
\E\Bigl[ \sup_{ s\le t} \bigl|\tilde{S}^1_s -
\tilde{S}^2_s\bigr|^2 \Bigr]^{1/2} \le
C_{\mmbox{Lem~\ref{lemSDELipschitzDependence}}} \rho_{\alpha\mbox{-\textup{\Holl}}}\bigl( \rrp
^1, \rrp^2
\bigr).
\]
\end{lemma}
\begin{pf}
Fix $R > 0$.
Let $\|\rrp^1\|_{\alpha\mbox{-\textup{\Holl}}}, \|\rrp^2\|_{\alpha\mbox{-\textup{\Holl}}} < R$.

From Lemma~\ref{lemTransformationRough} we know that
\[
\sup_{t,x} \bigl|\tilde{b}^1(t,x) - \tilde{b}^2(t,x)\bigr|
\le K_{\mmbox{Lem~\ref{lemTransformationRough}}}(R) \rho_{\alpha\mbox{-\textup{\Holl}}}\bigl(
\rrp^1,
\rrp^2 \bigr).
\]

Analogously, we get
\[
\sup_{t,x} \bigl|\tilde{a}^1(t,x) - \tilde{a}^2(t,x)\bigr|
\le L_2 \rho_{\alpha\mbox{-\textup{\Holl}}}\bigl( \rrp^1,
\rrp^2 \bigr)
\]
for a $L_2 = L_2(R)$.
\end{pf}
%
\begin{lemma}
\label{lemLipschitnessOfFlow}
Let $\alpha\in(0,1)$.
Let $\gamma> \frac{1}{\alpha} \ge1$, $k \in\{1, 2, \ldots\}$ and
assume that $V = (V_1,\ldots, V_d)$ is a collection of $\Lip^{\gamma
+k}$-vector fields on $\R^e$.
Write $n = (n_1,\ldots, n_e) \in\N^e$
and assume $|n|:= n_1 + \cdots+ n_e \le k$.

Then, for all $R > 0$ there exist $C = C(R, |V|_{\Lip^{\gamma+k}}), K
= K(R, |V|_{\Lip^{\gamma+k}}))$ such that if
$\x^1, \x^2 \in C^{ \alpha\mbox{-\textup{\Holl}}}( [0,t], G^{[p]}(\R^d) )$
with $\max_i \|\x^i\|_{\alpha\mbox{-\textup{\Holl}};[0,t]} \le R$, then
\begin{eqnarray*}
\sup_{y_0 \in\R^e} \bigl\llvert\partial_n \pi_{(V)}
\bigl(0,y_0; \x^1\bigr) - \partial_n
\pi_{(V)}\bigl(0,y_0; \x^2\bigr) \bigr
\rrvert_{\alpha\mbox{-\textup{\Holl}};[0,t]} &\le& C \rho_{\alpha\mbox{-\textup{\Holl}}}\bigl(\x^1,
\x^2\bigr),
\\
\sup_{y_0 \in\R^e} \bigl\llvert\partial_n \pi_{(V)}
\bigl(0,y_0; \x^1\bigr)^{-1} -
\partial_n \pi_{(V)}\bigl(0,y_0;
\x^2\bigr)^{-1} \bigr\rrvert_{\alpha\mbox{-\textup{\Holl}};[0,t]} &\le& C
\rho_{\alpha\mbox{-\textup{\Holl}}}\bigl(\x^1,\x^2\bigr),
\\
\sup_{y_0 \in\R^e} \bigl\llvert\partial_n \pi_{(V)}
\bigl(0,y_0; \x^1\bigr) \bigr\rrvert_{\alpha\mbox{-\textup{\Holl}}; [0,t]} &\le&
K,
\\
\sup_{y_0 \in\R^e} \bigl\llvert\partial_n \pi_{(V)}
\bigl(0,y_0; \x^1\bigr)^{-1} \bigr
\rrvert_{\alpha\mbox{-\textup{\Holl}}; [0,t]} &\le& K.
\end{eqnarray*}
\end{lemma}
\begin{pf}
The fact that $V \in\Lip^{\gamma+k}$
(instead of just $\Lip^{\gamma+k-1}$)
entails that the derivatives up to order $k$
are unique, nonexplosive solutions to RDEs with $\Lip^{\gamma
}_{\mathrm{loc}}$ vector fields; see Section 11 in~\cite{bibFrizVictoir}.
Localization (uniform for driving paths bounded in $\alpha$-H\"older norm)
then yields the desired results.
\end{pf}




\printaddresses

\end{document}